\documentclass[12pt]{article}
\usepackage{amsmath,amssymb,amsthm}
\usepackage{hyperref}
\usepackage{datetime}
\usepackage{units}
\usepackage{color}
\usepackage[T1]{fontenc}
\usepackage[utf8]{inputenc}
\usepackage{authblk}
\usepackage{bm,latexsym,mathrsfs,enumerate}
\usepackage{color}

\title{ A non-PSLQ route to BBP-type formulas\thanks{%
MSC 2010: 65B10, 11B99}}
\author[]{Kunle Adegoke\thanks{adegoke00@gmail.com}}
\affil{Department of Physics and Engineering Physics, \mbox{Obafemi Awolowo University, Ile-Ife, 220005 Nigeria}}

\theoremstyle{plain}
\numberwithin{equation}{section}

\begin{document}
\date{}
\maketitle
\begin{abstract}
\noindent BBP-type formulas are usually discovered experimentally, through computer searches. In this paper, however, starting with two simple generators, and hence without doing any computer searches, we derive a wide range of \mbox{BBP-type} formulas in general bases. Many previously discovered \mbox{BBP-type formulas} turn out to be particular cases of the formulas derived here.
\end{abstract}
\tableofcontents

\section{Introduction} 
The study of BBP-type (Bailey, Borwein and Plouffe, 1997) formulas has continued to attract attention, mainly because they facilitate digit extraction through a simple algorithm not requiring multiple-precision arithmetic~(Bailey, 2013). Experimentally, these formulas are usually discovered by using Bailey and Ferguson's PSLQ (Partial Sum of Squares -- Lower Quadrature) algorithm~(Ferguson \textit{et. al.}, 1999) or its variations. A downside is that PSLQ and other integer relation finding schemes typically do not suggest proofs~(Bailey, 2006). Formal proofs must be developed after the formulas have been discovered. There have been attempts in the past to give general formulas which include the proofs, as can be found, for example, in the following references:~(Bellard, 1997), (Broadhurst, 1998) and (Adamchik and Wagon, 1996). In this paper we give two identities which generate a wide range of BBP-type formulas in arbitrary bases. Many BBP-type formulas that are known in the literature turn out to be mere particular instances of the more general formulas presented here.

\section{Generators of BBP-type formulas}
Consider the Taylor series expansion
\begin{equation}\label{equ.mq2ge1q}
 - \ln (1 - z) = \sum_{k = 1}^\infty  {\frac{{z^k }}{k}}\,, 
\end{equation}
valid for $|z|\le 1$, $z\ne 1$.
Choosing $z=p\cos x+{\rm i}p\sin x$ in~\eqref{equ.mq2ge1q}, for real $p$ and $x$, allows one to write
\begin{align}
 - \ln (1 - z) &\equiv  - \ln \left[ {\sqrt {(1 - 2p\cos x + p^2 )} \exp \left( {{\rm i}\arctan \left( {\frac{{ - p\sin x}}{{1 - p\cos x}}} \right)} \right)} \right]\nonumber\\ 
 &=  - \frac{1}{2}\ln (1 - 2p\cos x + p^2 ) + {\rm i}\arctan \left( {\frac{{p\sin x}}{{1 - p\cos x}}} \right)\label{equ.zzd5y3t}\\
\intertext{and, using De Moiver theorem,}
\sum_{k = 1}^\infty  {\frac{{z^k }}{k}}  &\equiv \sum_{k = 1}^\infty  {\frac{{(p\cos x + {\rm i}p\sin x)^k }}{k} = \sum_{k = 1}^\infty  {\frac{{p^k \cos kx}}{k}} }  + {\rm i}\sum_{k = 1}^\infty  {\frac{{p^k \sin kx}}{k}}\,.\label{equ.vnk9uc0} 
\end{align}

Equating real and imaginary parts of~\eqref{equ.zzd5y3t} and~\eqref{equ.vnk9uc0} leads to the following identities:

\begin{equation}\label{equ.arctan}
\arctan \left( {\frac{{p\sin x}}{{1 - p\cos x}}} \right)=\sum_{k = 1}^\infty  {\frac{{p^k \sin {kx}}}{k}}
\end{equation}
and
\begin{equation}\label{equ.log}
- \frac{1}{2}\ln (1 - 2p\cos x + p^2 )=\sum_{k = 1}^\infty  {\frac{{p^k \cos {kx}}}{k}}\,.
\end{equation}

\bigskip

In the rest of this paper we demonstrate that careful choices of $p$ and $x$ in~\eqref{equ.arctan} and ~\eqref{equ.log} lead to interesting BBP-type series, for $|p|<1$. 

\section{Arctangent formulas}
\subsection{BBP-type formulas generated by $x=\pi/2$ in \mbox{identity}~\eqref{equ.arctan}}
The choice $x=\pi/2$ in~\eqref{equ.arctan} gives the identity
\begin{equation}\label{equ.vvky6mn}
\arctan p = \sum_{k = 1}^\infty  {\frac{{p^k \sin ({{k\pi } \mathord{\left/
 {\vphantom {{k\pi } 2}} \right.
 \kern-\nulldelimiterspace} 2})}}{k}}\,. 
\end{equation}
Since
\[
\sin\left(\frac{k\pi}{2}\right)=
\begin{cases}
1 & \text{if $k=1,5,9,13,17,\ldots$}\\
0 & \text{if $k=2,4,6,8,10,\ldots$}\\
-1 & \text{if $k=3,7,11,15,19,\ldots$}
\end{cases}\,,
\]
the identity~\eqref{equ.vvky6mn} can be written as

\[
\arctan p = \sum_{k = 0}^\infty  {p^{4k} \left[ {\frac{p}{{4k + 1}} - \frac{{p^3 }}{{4k + 3}}} \right]}\,. 
\]
Setting \mbox{$p=1/u$}, the above identity can be written as
\[
u^3 \arctan \frac{1}{u} = \sum_{k = 0}^\infty  {\frac{1}{{u^{4k} }}\left[ {\frac{{u^2 }}{{4k + 1}} - \frac{1}{{4k + 3}}} \right]}\,, 
\]
which is a BBP-type formula if $u^2$ is a positive integer. Thus, with $u^2=n$, we have:

\begin{equation}\label{equ.cm0jf0r}
\sqrt n \arctan\left( {\frac{1}{{\sqrt n }}} \right) = \frac{1}{n}\sum_{k = 0}^\infty  {\frac{1}{{(n^2 )^k }}\left[ {\frac{n}{{4k + 1}} - \frac{1}{{4k + 3}}} \right]},\quad n\in \mathbb{Z^+}\,.
\end{equation}
 
In the notation employed in the BBP Compendium~(Bailey, 2013),
\[
\sqrt n\arctan (\frac{1}{{\sqrt n }}) = \frac{1}{{n }}P(1,n^2 ,4,(n,0, - 1,0))\,.
\]
We note that the particular case $n=2$ is a \mbox{base-$4$} version of formula~(21) of the Compendium. To see this we write the base~$n^2$, length~$4$ formula~\eqref{equ.cm0jf0r} as a base~$n^4$, length~$8$ formula as follows:
\begin{equation*}
\begin{split}
&\sum_{k = 0}^\infty  {\frac{1}{{(n^2 )^k }}\left[ {\frac{n}{{4k + 1}} - \frac{1}{{4k + 3}}} \right]}  = \sum_{\mbox{$k$ even}} {\left( {} \right)}  + \sum_{\mbox{$k$ odd}} {\left( {} \right)}\\ 
&= \sum_{k = 0}^\infty  {\frac{1}{{(n^2 )^{2k} }}\left[ {\frac{n}{{4(2k) + 1}} - \frac{1}{{4(2k) + 3}}} \right]}  + \sum_{k = 0}^\infty  {\frac{1}{{(n^2 )^{2k + 1} }}\left[ {\frac{n}{{4(2k + 1) + 1}} - \frac{1}{{4(2k + 1) + 3}}} \right]}\\
& = \frac{1}{{n^2 }}\sum_{k = 0}^\infty  {\frac{1}{{(n^4 )^k }}\left[ {\frac{{n^3 }}{{8k + 1}} - \frac{{n^2 }}{{8k + 3}} + \frac{n}{{8k + 5}} - \frac{1}{{8k + 7}}} \right]}\,,  
\end{split}
\end{equation*}

so that in base~$n^4$, length~$8$ we have
\begin{equation}\label{equ.qpdxnm6}
\begin{split}
\sqrt n \arctan \left( {\frac{1}{{\sqrt n }}} \right)& = \frac{1}{{n^3 }}\sum_{k = 0}^\infty  {\frac{1}{{(n^4 )^k }}\left[ {\frac{{n^3 }}{{8k + 1}} - \frac{{n^2 }}{{8k + 3}} + \frac{n}{{8k + 5}} - \frac{1}{{8k + 7}}} \right]}\\ 
& = \frac{1}{{n^3 }}P(1,n^4 ,8,(n^3 ,0, - n^2 ,0,n, - 1))\,.
\end{split}
\end{equation}
The particular case $n=2$ in~\eqref{equ.qpdxnm6} recovers formula~(21) of the Compendium.

\bigskip

Similarly, the particular case $n=3$ in~\eqref{equ.cm0jf0r} is a \mbox{base-$9$} length~$4$ version of formula~(66) of the Compendium.

\bigskip

In general a formula with base $b$ and length $l$ can be rewritten as a formula with base $b^r$ and length $rl$ (Bailey, 2013).

\bigskip

Identity~\eqref{equ.vvky6mn} can also be written as
\[
\arctan p =p\sum_{k = 0}^\infty  {\frac{{( - p^2 )^k }}{{2k + 1}}}\,, 
\]
which gives the alternating base~$n$ version of~\eqref{equ.cm0jf0r} as
\[
\begin{split}
\sqrt n \arctan\left( {\frac{1}{{\sqrt n }}} \right) &= \sum_{k = 0}^\infty  {\frac{1}{{(-n)^k }}\left[ {\frac{1}{{2k + 1}}}\right]},\quad n\in \mathbb{Z^+}\\
&= P(1, - n,2,(1,0))\,.
\end{split}
\]

\subsection{BBP-type formulas generated by $x=\pi/3$ in \mbox{identity}~\eqref{equ.arctan}}
Putting $x=\pi/3$ in~\eqref{equ.arctan} gives
\[
\arctan \left( {\frac{{p\sqrt 3 }}{{2 - p}}} \right) = \sum_{k = 1}^\infty  {\frac{{p^k \sin ({{k\pi } \mathord{\left/
 {\vphantom {{k\pi } 3}} \right.
 \kern-\nulldelimiterspace} 3})}}{k}} \,.
\]

Noting that
\[
\sin\left(\frac{k\pi}{3}\right)=\frac{\sqrt 3}{2}
\begin{cases}
1 & $k=1,2,7,8,13,14,\ldots$\\
0 & $k=0,3,6,9,12,15,\ldots$\\
-1 & $k=4,5,10,11,16,17,\ldots$
\end{cases}\,,
\]
we have the following identity:
\[
\sqrt 3\arctan \left( {\frac{{p\sqrt 3 }}{{2 - p}}} \right) = \frac{{3 }}{2}\sum_{k = 0}^\infty  {( - 1)^k p^{3k} \left[ {\frac{p}{{3k + 1}} + \frac{{p^2 }}{{3k + 2}}} \right]}\,, 
\]
which is a BBP-formula if $p=1/\pm n$, $n$ a positive integer. The choice $p=1/n$ leads to

\begin{equation}\label{equ.kcbwvt9}
n^2 \sqrt 3\arctan \left( {\frac{{\sqrt 3 }}{{2n - 1}}} \right) = \frac{{3 }}{2}\sum_{k = 0}^\infty  {\frac{1}{{( - n^3 )^k }}\left[ {\frac{n}{{3k + 1}} + \frac{1}{{3k + 2}}} \right]}\,,\quad n\in \mathbb{Z^+}\,, 
\end{equation}
while the choice $p=-1/n$ gives
\begin{equation}\label{equ.n6pptf8}
n^2 \sqrt 3\arctan \left( {\frac{{\sqrt 3 }}{{2n + 1}}} \right) = \frac{{ 3 }}{2}\sum_{k = 0}^\infty  {\frac{1}{{( n^3 )^k }}\left[ {\frac{n}{{3k + 1}}- \frac{1}{{3k + 2}}} \right]}\,,\quad n\in \mathbb{Z^+}\,. 
\end{equation}

That is
\[
n^2 \sqrt 3\arctan \left( {\frac{{\sqrt 3 }}{{2n - 1}}} \right) = \frac{{ 3 }}{2}P(1, - n^3 ,3,(n,1,0))
\]
and
\[
n^2 \sqrt 3\arctan \left( {\frac{{\sqrt 3 }}{{2n + 1}}} \right) = \frac{{ 3 }}{2}P(1, n^3 ,3,(n,-1,0))\,.
\]

A particular case of~\eqref{equ.n6pptf8} is formula~(65) in the BBP Compendium, corresponding to $n=3$ here. $n=2$ in~\eqref{equ.kcbwvt9} also gives a formula that is equivalent to formula~(18) in the Compendium.

\subsection{BBP-type formulas generated by $x=\pi/4$ in \mbox{identity}~\eqref{equ.arctan}}
$x=\pi/4$ in~\eqref{equ.arctan} gives
\[
\arctan \left( {\frac{p}{{\sqrt 2  - p}}} \right) = \sum_{k = 1}^\infty  {\frac{{p^k \sin ({{k\pi } \mathord{\left/
 {\vphantom {{k\pi } 4}} \right.
 \kern-\nulldelimiterspace} 4})}}{k}}\,. 
\]
Observing that
\[
\sin\left(\frac{k\pi}{4}\right)=
\begin{cases}
1 & $k=2,10,18,26,34,\ldots$\\
1/\sqrt 2  & $k=1,3,9,11,17,19,\ldots$\\
0 & $k=0,4,8,12,16,20,\ldots$\\
-1/\sqrt 2 & $k=5,7,13,15,21,23,\ldots$\\
-1 & $k=6,14,22,30,38,46,\ldots$
\end{cases}\,,
\]
we obtain
\[
\begin{split}
\arctan \left( {\frac{p}{{\sqrt 2  - p}}} \right) = \sum_{k = 0}^\infty  {p^{8k} \left[ {\frac{p}{{\sqrt 2 }}\frac{1}{{8k + 1}} + \frac{{p^2 }}{{8k + 2}} + \frac{{p^3 }}{{\sqrt 2 }}\frac{1}{{8k + 3}}} \right.}\\
\left. { - \frac{{p^5 }}{{\sqrt 2 }}\frac{1}{{8k + 5}} - \frac{{p^6 }}{{8k + 6}} - \frac{{p^7 }}{{\sqrt 2 }}\frac{1}{{8k + 7}}} \right]\,.
\end{split}
\]
On setting $p=\sqrt 2 /u$, the above identity can be written as
\begin{equation}\label{equ.qdmzvaj}
\begin{split}
u^7 \arctan \left( {\frac{1}{{u - 1}}} \right) = \sum_{k = 0}^\infty  {\frac{1}{{({u \mathord{\left/
 {\vphantom {u {\sqrt 2 }}} \right.
 \kern-\nulldelimiterspace} {\sqrt 2 }})^{8k} }}\left[ {\frac{{u^6 }}{{8k + 1}} + \frac{{2u^5 }}{{8k + 2}} + \frac{{2u^4 }}{{8k + 3}}} \right.}\\
 - \left. {\frac{{4u^2 }}{{8k + 5}} - \frac{{8u}}{{8k + 6}} - \frac{8}{{8k + 7}}} \right]\,,
\end{split}
\end{equation}

while $p=-\sqrt 2 /u$ gives
\begin{equation}\label{equ.qoj5aln}
\begin{split}
u^7 \arctan \left( {\frac{1}{{u + 1}}} \right) = \sum_{k = 0}^\infty  {\frac{1}{{({u \mathord{\left/
 {\vphantom {u {\sqrt 2 }}} \right.
 \kern-\nulldelimiterspace} {\sqrt 2 }})^{8k} }}\left[ {\frac{{u^6 }}{{8k + 1}} - \frac{{2u^5 }}{{8k + 2}} + \frac{{2u^4 }}{{8k + 3}}} \right.}\\
 - \left. {\frac{{4u^2 }}{{8k + 5}} + \frac{{8u}}{{8k + 6}} - \frac{8}{{8k + 7}}} \right]\,.
\end{split}
\end{equation}

Identities~\eqref{equ.qdmzvaj} and~\eqref{equ.qoj5aln} are \mbox{BBP-type} series if $u$ is an even integer. Thus, setting $u=2n$ in both identities, we obtain the following \mbox{BBP-type} formulas:
\begin{equation}\label{equ.aeqgtu9}
\begin{split}
 n^7 \arctan \left( {\frac{1}{{2n - 1}}} \right) = \frac{1}{16}\sum_{k = 0}^\infty  {\frac{1}{{(16n^8 )^k }}\left[ {\frac{{8 n^6 }}{{8k + 1}} + \frac{{8 n^5 }}{{8k + 2}} + \frac{{4 n^4 }}{{8k + 3}}} \right.}\\
\qquad\left. { - \frac{{2 n^2 }}{{8k + 5}} - \frac{{2 n}}{{8k + 6}} - \frac{{1 }}{{8k + 7}}} \right],\quad n\in \mathbb{Z^+}
\end{split}
\end{equation}
and
\begin{equation}\label{equ.wneyxa4}
\begin{split}
n^7 \arctan \left( {\frac{1}{{2n + 1}}} \right) = \frac{1}{16}\sum_{k = 0}^\infty  {\frac{1}{{(16n^8 )^k }}\left[ {\frac{{8 n^6 }}{{8k + 1}} - \frac{{8 n^5 }}{{8k + 2}} + \frac{{4 n^4 }}{{8k + 3}}} \right.}\\
\qquad\left. { - \frac{{2 n^2 }}{{8k + 5}} + \frac{{2 n}}{{8k + 6}} - \frac{{1 }}{{8k + 7}}} \right],\quad n\in \mathbb{Z^+}\,.
\end{split}
\end{equation}

In the P-notation then,
\[
n^7 \arctan \left( {\frac{1}{{2n - 1}}} \right) = \frac{1}{{16}}P(1,16n^8 ,8,(8n^6 ,8n^5 ,4n^4 ,0, - 2n^2 , - 2n, - 1,0))
\]

and

\[
n^7 \arctan \left( {\frac{1}{{2n + 1}}} \right) = \frac{1}{{16}}P(1,16n^8 ,8,(8n^6 , - 8n^5 ,4n^4 ,0, - 2n^2 ,2n, - 1,0))\,.
\]
Formula~(15) of the Compendium is a particular case of~\eqref{equ.aeqgtu9}, with $n=1$.

\bigskip

Adding~\eqref{equ.qdmzvaj} and~\eqref{equ.qoj5aln}, we obtain
\[
u^7 \arctan \left( {\frac{{2u}}{{u^2  - 2}}} \right) = 2\sum_{k = 0}^\infty  {\frac{1}{{({u \mathord{\left/
 {\vphantom {u {\sqrt 2 }}} \right.
 \kern-\nulldelimiterspace} {\sqrt 2 }})^{8k} }}\left[ {\frac{{u^6 }}{{8k + 1}} + \frac{{2u^4 }}{{8k + 3}} - \frac{{4u^2 }}{{8k + 5}} - \frac{8}{{8k + 7}}} \right]}\,,
\]
which is a BBP-series only if $u^2$ is an even integer. Thus, setting $u^2=2n$ in the above identity, we obtain the BBP-type formula

\[
n^3 \sqrt {2n} \arctan \left( {\frac{{\sqrt {2n} }}{{n - 1}} } \right) = 2 \sum_{k = 0}^\infty  {\frac{1}{{n^{4k} }}\left[ {\frac{{n^3 }}{{8k + 1}} + \frac{{n^2 }}{{8k + 3}} - \frac{n}{{8k + 5}} - \frac{1}{{8k + 7}}} \right]}\,, 
\]
that is
\[
n^3 \sqrt {2n} \arctan \left( {\frac{{\sqrt {2n} }}{{n - 1}}} \right) = 2 P(1,n^4 ,8,(n^3 ,0,n^2 ,0, - n,0, - 1,0))\,.
\]
The particular case $n=2$ corresponds to formula~(8) in the Compendium. 

\bigskip

Subtracting~\eqref{equ.qoj5aln} from~\eqref{equ.qdmzvaj} gives a formula which is equivalent to~\eqref{equ.cm0jf0r} and therefore contains no new information.

\subsection{BBP-type formulas generated by $x=\pi/6$ in \mbox{identity}~\eqref{equ.arctan}}
With $x=\pi/6$ in~\eqref{equ.arctan}, we have
\[
\arctan \left( {\frac{{2p + p^2 \sqrt 3 }}{{4 - 3p^2 }}} \right) = \sum_{k = 1}^\infty  {\frac{{p^k \sin (k\pi /6)}}{k}}\,. 
\]
Noting that
\[
\sin\left(\frac{k\pi}{6}\right)=
\begin{cases}
1 & $k=3,15,27,39,51,\ldots$\\
\sqrt 3/2  & $k=2,4,14,16,26,28,\ldots$\\
1/2  & $k=1,5,13,17,25,29,37,41\ldots$\\
0 & $k=6,12,18,24,30,36\ldots$\\
-1/2  & $k=7,11,19,23,31,35,\ldots$\\
-\sqrt 3/2  & $k=8,10,20,22,32,34,\ldots$\\
-1 & $k=9,21,33,45,57,\ldots$
\end{cases}\,,
\]
we obtain
\[
\begin{split}
&\arctan \left( {\frac{{2p + p^2 \sqrt 3 }}{{4 - 3p^2 }}} \right)\\
&\qquad = \sum_{k = 0}^\infty  {( - p^6 )^k \left[ {\frac{1}{2}\frac{p}{{6k + 1}} + \frac{{\sqrt 3 }}{2}\frac{{p^2 }}{{6k + 2}} + \frac{{p^3 }}{{6k + 3}} + \frac{{\sqrt 3 }}{2}\frac{{p^4 }}{{6k + 4}} + \frac{1}{2}\frac{{p^5 }}{{6k + 5}}} \right]}\,.
\end{split}
\]
$p=\sqrt 3/u$ and $p=-\sqrt 3/u$ in the above identity yield the following series:
\[
\begin{split}
&u^5 \sqrt 3\arctan \left( {\frac{{\sqrt 3 }}{{2u - 3}}} \right)\\
&\qquad = \frac{{ 3 }}{2}\sum_{k = 0}^\infty  {\frac{1}{{({{ - u^6 } \mathord{\left/
 {\vphantom {{ - u^6 } {27}}} \right.
 \kern-\nulldelimiterspace} {27}})^k }}\left[ {\frac{{u^4 }}{{6k + 1}} + \frac{{3u^3 }}{{6k + 2}} + \frac{{6u^2 }}{{6k + 3}} + \frac{{9u}}{{6k + 4}} + \frac{9}{{6k + 5}}} \right]}
\end{split}
\]
and
\[
\begin{split}
&u^5 \sqrt 3\arctan \left( {\frac{{\sqrt 3 }}{{2u + 3}}} \right)\\
&\qquad= \frac{{ 3 }}{2}\sum_{k = 0}^\infty  {\frac{1}{{({{ - u^6 } \mathord{\left/
 {\vphantom {{ - u^6 } {27}}} \right.
 \kern-\nulldelimiterspace} {27}})^k }}\left[ {\frac{{u^4 }}{{6k + 1}} - \frac{{3u^3 }}{{6k + 2}} + \frac{{6u^2 }}{{6k + 3}} - \frac{{9u}}{{6k + 4}} + \frac{9}{{6k + 5}}} \right]}\,,
\end{split}
\]
which are BBP-type series if $u$ is a multiple of $3$. Thus, with $u=3n$, we obtain the following BBP-type series:
\begin{equation}\label{equ.hm0e9a2}
\begin{split}
& 27n^5 \sqrt 3\arctan \left( {\frac{1}{{\sqrt 3 }}\frac{1}{{2n - 1}}} \right)\\
&= \frac{{3 }}{2}\sum_{k = 0}^\infty  {\frac{1}{{( - 27n^6 )^k }}\left[ {\frac{{9n^4 }}{{6k + 1}} + \frac{{9n^3 }}{{6k + 2}} + \frac{{6n^2 }}{{6k + 3}} + \frac{{3n}}{{6k + 4}} + \frac{1}{{6k + 5}}} \right]}\\
&=\frac{{ 3 }}{2}P(1, - 27n^6 ,6,(9n^4 ,9n^3 ,6n^2 ,3n,1,0))
\end{split}
\end{equation}
and
\begin{equation}\label{equ.tggu3m5}
\begin{split}
& 27n^5 \sqrt 3\arctan \left( {\frac{1}{{\sqrt 3 }}\frac{1}{{2n + 1}}} \right)\\
&= \frac{{3 }}{2}\sum_{k = 0}^\infty  {\frac{1}{{( - 27n^6 )^k }}\left[ {\frac{{9n^4 }}{{6k + 1}} - \frac{{9n^3 }}{{6k + 2}} + \frac{{6n^2 }}{{6k + 3}} - \frac{{3n}}{{6k + 4}} + \frac{1}{{6k + 5}}} \right]}\\
&=\frac{{ 3 }}{2}P(1, - 27n^6 ,6,(9n^4 ,-9n^3 ,6n^2 ,-3n,1,0))\,.
\end{split}
\end{equation}

Formula~(66) of the Compendium is a particular case of formula~\eqref{equ.hm0e9a2}, corresponding to setting $n=1$.

\bigskip

Addition of~\eqref{equ.hm0e9a2} and~\eqref{equ.tggu3m5} gives the following BBP-type series
\begin{equation}
\begin{split}
n^2 \sqrt n \arctan \left( {\frac{{\sqrt n }}{{n - 1}}} \right) &= \sum_{k = 0}^\infty  {\frac{1}{{( - n^3 )^k }}\left[ {\frac{{n^2 }}{{6k + 1}} + \frac{{2n}}{{6k + 3}} + \frac{1}{{6k + 5}}} \right]}\\
&= P(1, - n^3 ,6,(n^2 ,0,2n,0,1,0)\,.
\end{split}
\end{equation}

Subtraction of~\eqref{equ.hm0e9a2} and~\eqref{equ.tggu3m5} yields~\eqref{equ.kcbwvt9} and therefore does not give new information. 
\section{Logarithm formulas}
Working in a similar fashion to that in the previous section, we present the following BBP-type formulas for logarithm. 
\subsection{BBP-type formulas generated by $x=\pi/2$ in \mbox{identity}~\eqref{equ.log}}
\begin{equation}\label{equ.miorglw}
\begin{split}
 \ln \left( {\frac{{n + 1}}{n}} \right) &=  \frac{1}{n^2}\sum_{k = 0}^\infty  {\frac{1}{{(n^2 )^k }}\left[ {\frac{n}{{2k + 1}} - \frac{1}{{2k + 2}}} \right]}\\
 &=  \frac{1}{n^2}P(1,n^2 ,2,(n,-1))\,.
\end{split}
\end{equation}
\begin{equation}\label{equ.dstx7ac}
\begin{split}
 \ln \left( {\frac{{n -1}}{n}} \right) &=  -\frac{1}{n^2}\sum_{k = 0}^\infty  {\frac{1}{{(n^2 )^k }}\left[ {\frac{n}{{2k + 1}} + \frac{1}{{2k + 2}}} \right]}\\
 &=  -\frac{1}{n^2}P(1,n^2 ,2,(n,1))\,.
\end{split}
\end{equation}
Addition of~\eqref{equ.miorglw} and~\eqref{equ.dstx7ac} gives
\begin{equation}\label{equ.wm2ipqk}
\begin{split}
\ln \left( {\frac{{n - 1}}{n}} \right) &=  - \frac{1}{n}\sum_{k = 0}^\infty  {\frac{1}{{n^k }}\left[ {\frac{1}{{k + 1}}} \right]}\\
& =  - \frac{1}{n}P(1,n,1,(1))\,.
\end{split}
\end{equation}

Fomula~(81) of the BBP Compendium is a particular case of~\eqref{equ.wm2ipqk}, with $n=10$.

\bigskip

Subtraction of~\eqref{equ.dstx7ac} from~\eqref{equ.miorglw} gives
\begin{equation}\label{equ.mxq900h}
\begin{split}
\sqrt n\;\ln \left( {\frac{{\sqrt n + 1}}{{\sqrt n - 1}}} \right) &= 2\sum_{k = 0}^\infty  {\frac{1}{{n^k }}\left[ {\frac{1}{{2k + 1}}} \right]}\\
 &= 2P(1,n ,2,(1,0))\,.
\end{split} 
\end{equation}
Formula~(6) of the Compendium is a particular case of~\eqref{equ.mxq900h}, with $n=4$, while Formula~(64) of the Compendium is another particular case, with $n=9$. $n=2$ in~\eqref{equ.mxq900h} gives a length~$2$ version of formula~(20) in the Compendium.
\subsection{BBP-type formulas generated by $x=\pi/3$ in \mbox{identity}~\eqref{equ.log}}
\[
\begin{split}
\ln \left( {\frac{{n^2  - n + 1}}{{n^2 }}} \right) &=  - \frac{1}{{n^3 }}\sum_{k = 0}^\infty  {\frac{1}{{( - n^3 )^k }}\left[ {\frac{{n^2 }}{{3k + 1}} - \frac{n}{{3k + 2}} - \frac{2}{{3k + 3}}} \right]}\\
& =  - \frac{1}{{n^3 }}P(1, - n^3 ,3,(n^2 , - n, - 2))\,. 
\end{split}
\]
\[
\begin{split}
\ln \left( {\frac{{n^2  + n + 1}}{{n^2 }}} \right) &= \frac{1}{{n^3 }}\sum_{k = 0}^\infty  {\frac{1}{{(n^3 )^k }}\left[ {\frac{{n^2 }}{{3k + 1}} + \frac{n}{{3k + 2}} - \frac{2}{{3k + 3}}} \right]}\\
& = \frac{1}{{n^3 }}P(1,n^3 ,3,(n^2 ,n, - 2))\,. 
\end{split}
\]

\subsection{BBP-type formulas generated by $x=\pi/4$ in \mbox{identity}~\eqref{equ.log}}

\[
\begin{split}
\ln \left( {\frac{{2n^2  - 2n + 1}}{{2n^2 }}} \right) &=  - \frac{1}{{2n^4 }}\sum_{k = 0}^\infty  {\frac{1}{{( - 4n^4 )^k }}\left[ {\frac{{2n^3 }}{{4k + 1}} - \frac{n}{{4k + 3}} - \frac{1}{{4k + 4}}} \right]}\\
&=  - \frac{1}{{2n^4 }}P(1, - 4n^4 ,4,(2n^3 ,0, - n, - 1))\,.
\end{split}
\]

\[
\begin{split}
\ln \left( {\frac{{2n^2  + 2n + 1}}{{2n^2 }}} \right) &= \frac{1}{{2n^4 }}\sum_{k = 0}^\infty  {\frac{1}{{( - 4n^4 )^k }}\left[ {\frac{{2n^3 }}{{4k + 1}} - \frac{n}{{4k + 3}} + \frac{1}{{4k + 4}}} \right]}\\ 
&= \frac{1}{{2n^4 }}P(1, - 4n^4 ,4,(2n^3 ,0, - n,1))\,.
\end{split}
\]

\begin{equation}\label{equ.tu1nql9}
\begin{split}
n\frac{\sqrt n}{\sqrt 2} \ln \left( {\frac{{n + \sqrt 2\sqrt n  + 1}}{{n - \sqrt 2\sqrt n  + 1}}} \right)& = 2\sum_{k = 0}^\infty  {\frac{1}{{( - n^2 )^k }}\left[ {\frac{n}{{4k + 1}} - \frac{1}{{4k + 3}}} \right]}\\
&= 2P(1, - n^2 ,4,(n,0,-1,0))\,.
\end{split}
\end{equation}
Note that $n=2$ in~\eqref{equ.tu1nql9} gives a binary BBP-type formula for $\log 5$.

\subsection{BBP-type formulas generated by $x=\pi/6$ in \mbox{identity}~\eqref{equ.log}}
\[
\begin{split}
\ln \left( {\frac{{3n^2  - 3n + 1}}{{3n^2 }}} \right) &=  - \frac{1}{{27n^6 }}\sum_{k = 0}^\infty  {\frac{1}{{( - 27n^6 )^k }}\left[ {\frac{{27n^5 }}{{6k + 1}} + \frac{{9n^4 }}{{6k + 2}} - \frac{{3n^2 }}{{6k + 4}} - \frac{{3n}}{{6k + 5}} - \frac{2}{{6k + 6}}} \right]}\\ 
&=  - \frac{1}{{27n^6 }}P(1, - 27n^6 ,6,(27n^5 ,9n^4 ,0,-3n^2 ,-3n,-2))\,.
\end{split}
\]
\[
\begin{split}
\ln \left( {\frac{{3n^2  + 3n + 1}}{{3n^2 }}} \right) &= \frac{1}{{27n^6 }}\sum_{k = 0}^\infty  {\frac{1}{{( - 27n^6 )^k }}\left[ {\frac{{27n^5 }}{{6k + 1}} - \frac{{9n^4 }}{{6k + 2}} + \frac{{3n^2 }}{{6k + 4}} - \frac{{3n}}{{6k + 5}} + \frac{2}{{6k + 6}}} \right]}\\
& = \frac{1}{{27n^6 }}P(1, - 27n^6 ,6,(27n^5 , - 9n^4 ,0,3n^2 , - 3n,2))\,.
\end{split}
\]
\begin{equation}\label{equ.ubf2fyh}
\begin{split}
n^2\frac{\sqrt n}{\sqrt 3} \ln \left( {\frac{{n + \sqrt 3\sqrt n  + 1}}{{n - \sqrt 3\sqrt n  + 1}}} \right)& = 2\sum_{k = 0}^\infty  {\frac{1}{{( - n^3 )^k }}\left[ {\frac{n^2}{{6k + 1}} - \frac{1}{{6k + 5}}} \right]}\\
&= 2P(1, - n^3 ,6,(n^2,0,0,0,-1,0))\,.
\end{split}
\end{equation}

The particular case $n=3$ of~\eqref{equ.ubf2fyh} is equivalent to the series obtained for $\log 7$ in reference~(Adamchik and Wagon, 1996).




\section{Summary}
Starting with two simple generators, we have derived a wide range of \mbox{BBP-type} formulas in general bases, namely, 
\[
\begin{split}
\sqrt n\arctan \left(\frac{1}{{\sqrt n }}\right) &= \frac{1}{{n }}P(1,n^2 ,4,(n,0, - 1,0))\\
n^2\sqrt 3\arctan \left( {\frac{{\sqrt 3 }}{{2n - 1}}} \right) &= \frac{{ 3 }}{2}P(1, - n^3 ,3,(n,1,0))\\
n^2\sqrt 3\arctan \left( {\frac{{\sqrt 3 }}{{2n + 1}}} \right) &= \frac{{ 3 }}{2}P(1, n^3 ,3,(n,-1,0))\,\\
n^7 \arctan \left( {\frac{1}{{2n - 1}}} \right) &= \frac{1}{{16}}P(1,16n^8 ,8,(8n^6 ,8n^5 ,4n^4 ,0, - 2n^2 , - 2n, - 1,0))\\
n^7 \arctan \left( {\frac{1}{{2n + 1}}} \right) &= \frac{1}{{16}}P(1,16n^8 ,8,(8n^6 , - 8n^5 ,4n^4 ,0, - 2n^2 ,2n, - 1,0))\\
n^3 \sqrt {2n} \arctan \left( {\frac{{\sqrt {2n} }}{{n - 1}} } \right) &= 2 P(1,n^4 ,8,(n^3 ,0,n^2 ,0, - n,0, - 1,0))\\
9n^5 \sqrt 3\arctan \left( {\frac{1}{{\sqrt 3 }}\frac{1}{{2n - 1}}} \right) &= \frac{{ 1 }}{2}P(1, - 27n^6 ,6,(9n^4 ,9n^3 ,6n^2 ,3n,1,0))\\
9n^5 \sqrt 3\arctan \left( {\frac{1}{{\sqrt 3 }}\frac{1}{{2n + 1}}} \right) &=\frac{{ 1 }}{2}P(1, - 27n^6 ,6,(9n^4 ,-9n^3 ,6n^2 ,-3n,1,0))\\
n^2 \sqrt n \arctan \left( {\frac{{\sqrt n }}{{n - 1}}} \right) &= P(1, - n^3 ,6,(n^2 ,0,2n,0,1,0)
\end{split}
\]
\[
\begin{split}
\ln \left( {\frac{{n + 1}}{n}} \right)&=  \frac{1}{n^2}P(1,n^2 ,2,(n,-1))\\
\ln \left( {\frac{{n -1}}{n}} \right) &=-\frac{1}{n^2}P(1,n^2 ,2,(n,1))\\ 
\sqrt n\ln \left( {\frac{{\sqrt n + 1}}{{\sqrt n - 1}}} \right) &= 2P(1,n ,2,(1,0))\\
\ln \left( {\frac{{n^2  + n + 1}}{{n^2 }}} \right)  &= \frac{1}{{n^3 }}P(1,n^3 ,3,(n^2 ,n, - 2))\\
\ln \left( {\frac{{n^2  - n + 1}}{{n^2 }}} \right) &=- \frac{1}{{n^3 }}P(1, - n^3 ,3,(n^2 , - n, - 2))\\ 
\ln \left( {\frac{{2n^2  - 2n + 1}}{{2n^2 }}} \right) &=  - \frac{1}{{2n^4 }}P(1, - 4n^4 ,4,(2n^3 ,0, - n, - 1))\\
\ln \left( {\frac{{2n^2  + 2n + 1}}{{2n^2 }}} \right) &=\frac{1}{{2n^4 }}P(1, - 4n^4 ,4,(2n^3 ,0, - n,1))\\
n\frac{\sqrt n}{\sqrt 2} \ln \left( {\frac{{n + \sqrt 2\sqrt n  + 1}}{{n - \sqrt 2\sqrt n  + 1}}} \right)&=2 P(1, - n^2 ,4,(n,0,-1,0))\\
\ln \left( {\frac{{3n^2  \pm 3n + 1}}{{3n^2 }}} \right) &=  \pm \frac{1}{{27n^6 }}P(1, - 27n^6 ,6,(27n^5 ,\mp 9n^4 ,0,\pm 3n^2 ,-3n,\pm 2))\\
n^2\frac{\sqrt n}{\sqrt 3} \ln \left( {\frac{{n + \sqrt 3\sqrt n  + 1}}{{n - \sqrt 3\sqrt n  + 1}}} \right)&= 2P(1, - n^3 ,6,(n^2,0,0,0,-1,0))
\end{split}
\]
Many previously discovered BBP-type formulas turn out to be particular cases of the above formulas.

\section*{References}

D.~H. Bailey, P.~B. Borwein and S.~Plouffe (1997), On the rapid computation of various polylogarithmic constants,
\emph{Mathematics of Computation}, 66 (218):903--913.\par
\bigskip
D.~H.~Bailey (2013), A compendium of BBP-type formulas for mathematical constants, [Online] Available:\\
\verb+ http://crd.lbl.gov/~dhbailey/dhbpapers/bbp-formulas.pdf+\par
\bigskip
H.~R.~P.~Ferguson, D.~H.~Bailey and S.~Arno (1999), Analysis of PSLQ, an integer relation finding algorithm,
\emph{Mathematics of Computation}, 68:351--369.\par
\bigskip
D.~H.~Bailey (2006), Algorithms for Experimental Mathematics I, [Online] Available:\\
\verb+ http://crd.lbl.gov/~dhbailey/dhbpapers/+\par
\bigskip
F.~Bellard (1997), A new formula to compute the $n^{th}$ binary digit of $pi$, [Online] Available:\\
\verb+ http://www-stud.enst.fr/~bellard/pi/pi_bin/pi_bin.html+\par
\bigskip
D.~J.~Broadhurst (1998), Polylogarithmic ladders, hypergeometric series and the ten millionth
  digits of $\zeta(3)$ and $\zeta(5)$.\\
\emph{arXiv:math/9803067v1 [math.CA]}\par
\bigskip
V.~Adamchik and S.~Wagon (1996), Pi: A 2000-year-old search changes direction,
\emph{Mathematica in Science and Education}, 5:11--19.\par

\end{document}